\documentclass[12pt, twoside]{article}
 \usepackage{amsmath,amsfonts,theorem}
 \numberwithin{equation}{section}
 
 \newtheorem{prop}{Proposition}[section]

 \theorembodyfont{\rmfamily}
 \newtheorem{dfn}{Definition}[section]

 \newcommand{\qed}{\ifhmode\unskip\nobreak\fi\quad\ensuremath\square}

 \newcommand{\ie}{i.\;e.\ }
 \newcommand{\wtd}{\widetilde}
 \newcommand{\Int}{\mathrm{Int}}
 \newcommand{\res}{{\mathrm{res}}}
 \newcommand{\bires}{{\mathrm{bi\,res}}}
 \newcommand{\bydef}{\,\stackrel{\mathrm{def}}{=}\,}
 
 \newcommand{\p}{\partial}
 \newcommand{\ov}{\overline}

 \newcommand{\sC}{\mathcal C}
 \newcommand{\sF}{\mathcal F}

 \newcommand{\sM}{\mathcal M}
 \newcommand{\Oh}{\mathcal O}
 \newcommand{\sP}{\mathcal P}
 
 \newcommand{\om}{\omega}
 \newcommand{\ro}{\varrho}
 \newcommand{\Ga}{\Gamma}
 \newcommand{\Om}{\Omega}
 \newcommand{\Si}{\Sigma}
 \newcommand{\PP}{\mathbb P}
 \newcommand{\C}{\mathbb C}

 \newcommand{\ad}{\operatorname{\mathrm{ad}}}
 \newcommand{\MSS}{\mathcal M^{\mathrm{ss}}}
 \newcommand{\diag}{\operatorname{\mathrm{diag}}}
 
 \newcommand{\CLRep}{\mathfrak{Rep}}
 \newcommand{\SL}{\operatorname{\text{\sl SL\/}}}
 \newcommand{\SU}{\operatorname{\text{\sl SU\/}}}
 \newcommand{\Map}{\text{\sl Map\/}}
 
 \newcommand{\Pic}{\operatorname{Pic}}

\begin{document}

 \title{{\huge\bf Hitchin systems on ll-curves}}
 \markboth{\hfill{\footnotesize A. N. Tyurin}\quad}
          {\quad{\footnotesize Hitchin systems on ll-curves}\hfill}

 \author{\fbox{\sc Andrei Tyurin}}
 \date{{\small 24.02.1942 -- 27.10.2002}}
 \maketitle

 \bigskip
 \begin{flushright}
 {\em To Prof. Friedrich Hirzebruch on his 75\,th birthday}
 \end{flushright}

 \bigskip

\section{Introduction}

Recall (see, for example \cite{T1} and \cite{T2}) that the
Deligne-Mumford compactification $\ov{\sM_{g}}$ of the moduli space
$\sM_{g}$ of smooth curves of genus $g$ containes a finite
configuration of points $\sP \subset \ov{\sM_{g}}$ corresponding to
the large limit curves and enumerated by 3-valent graphs. Given a
3-valent graph $\Ga$, we write:
\begin{itemize}
\item
$V(\Ga) = \{ v_i \}$ for the set of its vertices;
\item
$E(\Ga) = \{ e \}$ for the set of its edges;
\item
$\vec E(\Ga) = \{ \vec e \}$ for the set of all orientation choices
for all the edges;
\item
$v_s(\vec e)$, $v_t(\vec e)$ for the source and target vertices of
an oriented edge $\vec e\in\vec E$;
\item
$S(v)\subset E(\Ga)$ for the star of a vertex $v \in V(\Ga)$,
i.\;e. for the set of all (non oriented) edges incident to $v$;
\end{itemize}
Topologically, each 3-valent graph $\Ga$ is equivalent to some
3-dimensional handlebody $H_{\Ga}$ whose boundary $\p H_{\Ga} =
\Si_{\Ga}$ is a smooth compact Riemann surface obtained from $\Ga$
by pumping its edges to tubes and vertices to {\it trinions\/},
that is 2-spheres with 3 holes. It is easy to see that the
cardinalities $|E(\Ga)| = 3g-3$, $|V(\Ga)|=2g-2$, where $g$ is the
genus of $\Si_{\Ga}$. By the construction, $\Si_{\Ga}$ is equipped
with a {\it trinion decomposition\/}, \ie it is sliced into `pairs
of pants' obtained by cutting some meridian $C_e$ out of each tube
$e\in E(\Ga)$. Contracting these $3g-3$ meridians to  points, we
obtain a connected reducible algebraic curve $P_\Ga$ glued from
$2g-2$ Riemann spheres $P_{v}=\C\PP^{1}$, $v \in V(\Ga)$, along
triples of points $(p_{e_1},p_{e_2},p_{e_3})\in P_{v}$
corresponding to the edges of the star $S(v)=\{e_1,e_2,e_3\}$. This
reducible algebraic curve $P_\Ga$ is Deligne - Mumford
stable and has the arithmetical genus $g$.

The restriction of the canonical sheaf $\Oh(K_{\Ga})$ to each component
$P_{v}$ coincides with the  sheaf of meromorphic differentials
$\om$ with  simple poles at $(p_{e_1},p_{e_2},p_{e_3})\in
P_{v}$:
$$\Oh_{P_\Ga}(K_\Ga)|_{P_v}
  =\Oh_{P_v}(K_{P_v}+p_{e_1}+p_{e_2}+p_{e_3})
  =\Oh_{P_v}(1)\;.
$$
Thus, a global holomorphic section $s\in H^0(P_\Ga,\Oh(K_{\Ga}))$ is a
collection of meromorphic differentials $\{\om_{v}\}$ on $P_{v}$
whose residues at the poles satisfy the following compatibility relations:
\begin{equation}\label{rescond1}
\begin{gathered}
\res_{p_{e}}\om_{v_s(e)}+\res_{p_{e}}\om_{v_t(e)}=0
\,,\;\forall\;e\in E(\Ga)\\
\sum_{e\in S(v)}\res_{p_e}\om(v)=0
\,,\;\forall\;v\in V(\Ga)
\end{gathered}
\end{equation}
These linear constraints lead to the right value for $\dim
H^0(P_\Ga,\Oh(K_{\Ga}))=g$, the dimension of the space of global
holomorphic differentials. This agrees with the Deligne -  Mumford
prediction that the vector bundle $\pi:V_1\to\ov{\sM_g}$, whose
fiber at a smooth $C\in\sM_g$ is $H^0(C,\Oh(K_C))$, can be continued
into each ll-curve point as a vector bundle of rank $g$. However,
the properties of the complete canonical linear system $|K_\Ga|$
depend on the topology of $\Ga$ (more precisely, on the thickness
of the graph, see \cite{A}).

The double canonical system $|2K_{\Ga}|$ of an ll-curve is much
more regular. Namely, a holomorphic quadratic differential $\Om\in
H^0(P_\Ga,\Oh(2K_{\Ga})$ is given by a collection of meromorphic
quadratic differentials $\{\Om_{v}\}$ on the components $P_{v}$
such that their bi-residues satisfy the following condition:
\begin{equation}\label{rescond2}
\bires_{p_{e}}\om_{v_1}=\bires_{p_{e}}\om_{v_2}\;,
\quad\forall\;e=(v_1,v_2)\in E(\Ga)
\end{equation}
(see \cite{T4},\cite{T5}). The system of nodes $\{p_e\}\in P_\Ga$,
$e\in E(\Ga)$, gives an identification
\begin{equation}\label{2Kdec}
H^0(P_\Ga,\Oh(2K_\Ga))^* =\C^{E(\Ga)}
\end{equation}
under which the basic linear form $H_e$, which corresponds to $e\in
E(\Ga)$, takes  quadratic differential $\Om$ to
$H_e(\Om)=\bires_{p_e}\Om\;$. So, the moduli space $\ov{\sM_{g}}$
has a smooth orbifold structure at each ll-point $P_\Ga\in\ov{\sM_{g}}$
and the fiber of the tangent bundle $T_{P_{\Ga}} \ov{\sM_{g}}$
at $P_\Ga$ equals
\begin{equation}\label{tb}
T_{P_{\Ga}} \ov{\sM_{g}}=H^{0}(P_{\Ga},\Oh_{P_\Ga}(2K_{P_\Ga})^*
                        =\C^{E(\Ga)}
\end{equation}
Moreover, each basic tangent direction $\C\cdot H_e$, $e\in
E(\Ga)$, can be integrated to some rational curve
$\sC_{e\subset\Ga}=\psi_{e\subset\Ga}(\PP^{1})\subset \ov{\sM_g} $
(see (3.1) -- (3.12) in \cite{T2}).

So, the Deligne - Mumford  compactification $\ov{\sM_{g}}$ contains
the configuration of points $\sP=\{P_\Ga\}$ parameterized by
3-valent graphs $\Ga$ of genus $g$ and the configuration of
rational curves\footnote{we consider $\sC=\cup\sC_{e\subset\Ga}$ as
one reducible curve}
$\sC=\bigcup\limits_{e\subset\Ga}\sC_{e\subset\Ga}$ corresponding
to flags $e\subset\Ga$ in such a way that three flags from the same
nest (see (2.17), (2.20) in \cite{T2}) lead to the same curve
$\sC_{e\subset\Ga}$. In the previous paper \cite{T2} we described
the complex gauge theory on ll-curves $\sP\subset \ov{\sM_g}$. Here
we investigate the Hitchin systems on ll-curves.

\section{Framed vector bundles on ll-curves}

We write $\MSS(P_\Ga)$ for the moduli space of topologically
trivial semistable holomorphic vector bundles of rank 2 on $P_\Ga$.
By the definition, $\MSS(P_\Ga)$ parameterizes rank 2 vector
bundles $E$ that have no positive line sub bundles $L\subset E$ and
are restricted to trivial holomorphic vector bundles $E|_{P_v}$
over each sphere $P_v$, $v\in V(\Ga)$. Any topologically trivial
rank 2 vector bundle $E$ on $P_\Ga$ can be {\it framed\/} by fixing
some trivialization
$$E|_{P_v}=P_v\times\C^2=P_v\times V_0
$$
over each component of $P_v\subset P_\Ga$ (here $V_0$ is some fixed
2-dimensional vector space). We write $\sF(P_\Ga)$ for the
space of framed topologically trivial vector bundles on $P_\Ga$.
The framing forgetful map
\begin{equation}\label{fgmap}
f\,:\;\sF(\Ga)\to\MSS(P_\Ga)
\end{equation}
gives a principal bundle with the `complex gauge' group
\begin{equation}\label{cggdef}
\Map(V(\Ga),\SL(2,\C))\;,
\end{equation}
where $\SL(2,\C)=\SL(V_0)$ consists of linear
automorphisms of $V_0$ with determinant 1.

Each framing of $E$ can be considered as a function
\begin{equation}\label{afun}
a:\vec E(\Ga)\to\SL(2,\C)\text{ \small such that }
a(\overleftarrow{e})=a(\overrightarrow e)^{-1}
\end{equation}
where $\overleftarrow{e}$, $\overrightarrow e$ are
two opposite orientation choices for  edge $e\in E(\Ga)$.
In these terms, complex gauge transformation
$g\in \Map(V(\Ga),\SL(2,\C))$ acts on a framing
$a$ by the formula
\begin{equation}\label{gta}
g(a)(\vec e)=g(v_s)\cdot a(\vec e)\cdot g(v_t)^{-1}.
\end{equation}
Factorizing \eqref{fgmap} through this action, we get
\begin{equation}\label{quot1}
\MSS(P_\Ga)=\sF(\Ga)/\Map(V(\Ga),\SL(2, \C)).
\end{equation}
At the same time, \eqref{afun} and \eqref{gta} show that a
function $a$ from \eqref{afun} is nothing but a flat
$\SL(2,\C)$-connection on the graph $\Ga$. Moreover, the reframing
group \eqref{cggdef} is really the complex gauge group
acting on these connections. Thus, the quotient
\eqref{quot1} coincides with the space of the gauge orbits of flat
connections, which, in its own turn, is the same as the space
$\CLRep(\pi_1(\Ga),\SL(2, \C))$, of
equivalence classes of representations of the fundamental group
$\pi_1(\Ga)$. In other words,
\begin{equation}\label{mssasclrep}
  \MSS(P_\Ga)=\CLRep(\pi_1(\Ga),\SL(2, \C))\;.
\end{equation}

Let us write $F_g$ for the free group with $g$ generators. Then
there is an isomorphism of fundamental groups
$$\pi_1(\Ga)=\pi_1(H_\Ga) = F_g
$$
where $H_\Ga$ is the handlebody whose boundary
$\p H_{\Ga}=\Si_{\Ga}$ is the pumped graph. So,
the space of representation classes \eqref{mssasclrep} is
described as the factor of the direct product of $g$ copies of
$\SL(2,\C)$ through diagonal action of $\SL(2,\C)$ by conjugations:
\begin{equation}\label{shottkydef}
 \CLRep(\pi_1(\Ga),\SL(2,\C))=(\SL(2,\C))^g/\ad_{\diag}\SL(2,\C)
\end{equation}
The right hand side here is called {\it the Schottki space\/}
of genus $g$ and denoted by $S_g$. So, $\MSS(P_\Ga) = S_g$.

On the other side, we have the canonical surjection
\begin{equation}\label{rdef}
r\,:\;\pi_1(\Si_\Ga)\to\pi_1(H_\Ga)=\pi_1(\Ga)\,.
\end{equation}
Let us choose some standard presentation
\begin{equation}\label{genchoice}
  \pi_1(\Si_\Ga)=
  \bigl\langle\,a_1,\,\ldots\,,\,a_g,\,b_1,\,\ldots\,,\,b_g\,|\;
             \prod_{i=1}^g[a_i,b_i]=1\bigr\rangle
\end{equation}
such that
\begin{equation}\label{afix}
\begin{aligned}
\ker(r)&=\langle a_1,\,\ldots\,,\,a_g\rangle\simeq F_g\\
\pi_1(\Ga)=\pi_1(H_\Ga)
  &=\langle r(b_1),\,\ldots\,,\,r(b_g)\rangle\simeq F_g
\end{aligned}
\end{equation}
Then there are two Schottki subspaces
$S_g^a,S_g^b\subset\CLRep(\pi_1(\Si_\Ga),\SL(2, \C))$
associated with the presentation \eqref{genchoice}:
\begin{equation}\label{shpair}
\begin{aligned}
S_g^a&=\{\ro\in\CLRep(\pi_1(\Si_\Ga),\SL(2,\C))\,|\;
         \ro(a_i)=1\text{ {\small for} }1\le i\le g\}\\
S_g^b&=\{\ro\in\CLRep(\pi_1(\Si_\Ga),\SL(2,\C))\,|\;
         \ro(b_i)=1\text{ {\small for} }1\le i\le g\}
\end{aligned}
\end{equation}
Further, each oriented edge $\vec e\in\vec E(\Ga)$ gives a homotopy
class of oriented 1-cycle $\vec C_e$ presented by some (oriented)
meridian rounding about the tube corresponding to this edge.
Joining these meridians with some fixed base point $p\in
P_{v_0}\subset\Si_\Ga$, which lives on the pumped vertex $v_0\in
V(\Ga)$ say, we get elements $[C_{\vec e}]$ of the pointed
fundamental group $\pi_1(\Si_\Ga)_p$. These $3g-3$ cycles $[C_{\vec
e}]$ obviously lay in the kernel $\ker(r)$ of the projection
\eqref{rdef}.

Let us realize $\pi_1(\Ga)_{v_0}$ as the group of
oriented cycles $(\vec e_1,\,\ldots\,,\,\vec e_d)$, $v_t(\vec
e_i)=v_s(\vec e_{i+1})$ based at $v_0=v_s(\vec e_1)=v_t(\vec
e_d)$. Then there is a map
$$\Int:\pi_1(\Ga)_{v_0}\to\pi_1(\Si_\Ga)_p
$$
defined as $\Int(\vec e_1,\,\ldots\,,\,\vec e_d)\bydef
  [C_{\vec e_1}]\circ[C_{\vec e_2}]\circ\,\cdots\,\circ
  [C_{\vec e_d}]$.
Clearly, $\Int(\pi_1(\Ga))=\ker(r)$ and the exact triple of
fundamental groups:
$$1\to\pi_1(\Ga)\stackrel{\Int}{\to}\pi_1(\Si_\Ga)\stackrel{r}{\to}
    \pi_1(\Ga)\to 1.
$$
induces a chain of maps between the corresponding spaces of
equivalence classes of $\SL(2,\C)$-representations:
\begin{equation}\label{reptrip}
\CLRep(\pi_1(\Ga))\stackrel{r^*}{\to}
\CLRep(\pi_1(\Si_\Ga))\stackrel{\Int^*}{\to}
\CLRep(\pi_1(\Ga))
\end{equation}
which is `exact' in the sense that each composition $r^*\circ\Int^*$
is the constant map into the identity. The middle space
$\CLRep(\pi_1(\Si_\Ga),\SL(2,\C))$ in \eqref{reptrip} parameterizes
topologically trivial {\it flat\/} holomorphic bundles on $P_\Ga$,
that is the pairs $(E,h)$, where $E\in\MSS(P_\Ga)$ and $h$ is a
holomorphic flat $\SL(2,\C)$-connection on $E$ (see formula (4.16)
and Proposition 4.1 from \cite{T2}). Moreover, the following
identification list can be checked at once:

\begin{prop}\ \par\vspace{-1ex}
  \begin{enumerate}
    \item
    The coincidence $\CLRep(\pi_1(\Ga))=\MSS(P_\Ga)$ described in
    \eqref{mssasclrep} identifies the second map of
    \eqref{reptrip}\,:\quad
    $\CLRep(\pi_1(\Si_\Ga))\stackrel{\Int^*}{\to}\CLRep(\pi_1(\Ga))$\quad
    with the forgetful map
    \begin{equation}\label{hvb}
       \CLRep(\pi_1(\Si_\Ga))\stackrel{f}{\to}\MSS(P_\Ga)\;,
    \end{equation}
    which takes $(E,h)\longmapsto E$. In particular, the fibers of
    $f=\Int^*$ have the structure of affine spaces associated with
    the vector spaces of Higgs fields (see for example \upshape{\cite{T4},
    \cite{T5}}).
    \item
    The first map of \eqref{reptrip}\,:\quad
    $\MSS(P_\Ga)=\CLRep(\pi_1(\Ga))\stackrel{r^*}{\to}
      \CLRep(\pi_1(\Si_\Ga))$\quad
    gives a section for the affine bundle \eqref{hvb}. So,
    each fiber of $f=\Int^*$ is naturally identified with the vector
    space of the Higgs fields on the corresponding vector bundle $E$.
  \end{enumerate}
\end{prop}

\noindent
There is quite simple geometrical reason for why
$\CLRep(\pi_1(\Si_\Ga))\to\MSS(P_\Ga)$ has a vector bundle
structure: the boundary $\ov{\MSS(P_\Ga)}\setminus\MSS(P_\Ga)$ of
the compactifified moduli space contains an effective theta divisor
(see \cite{T1}); so, the obstruction for lifting of affine
structure to  vector one vanishes over $\MSS(P_\Ga)$.

\section{Framed Higgs fields and spectral curves}

Let us write $\wtd E_a$ for a framed vector bundle obtained from
$E$ by a framing function $a$ as in \eqref{afun}. Then a Higgs
field $\phi:\wtd E_a\to\wtd E_a\otimes K_\Ga$ on  $\wtd
E_a$ is given by a collection of traceless
$2\times2$-matrices $\om(v)|$, $v\in V(\Ga)$,
whose entries $\om_{ij}(v)$ are meromorphic
differentials on $P_v$ with poles at $p_{e_1}$, $p_{e_2}$,
$p_{e_3}$, $\{e_1,e_2,e_3\}=S(v)$. Their residue matrices
$\res_{p_e}(\om(v))=||\res_{p_e}\om_{ij}(v)||$ satisfy the relations:
\begin{equation}\label{resmatcond}
\begin{gathered}
  \res_{p_e}\om(v_s(\vec e))+
   a(\vec e)\circ\res_{p_e}\om(v_t(\vec e))\circ a(\vec e)^{-1}=0
  \,,\;\forall\;\vec e\in\vec E(\Ga)\\
  \sum_{e\in S(v)}\res_{p_e}\om(v)=0
  \,,\;\forall\;v\in V(\Ga)\,.
\end{gathered}
\end{equation}
These constraints are gauge invariant\footnote{with respect to the
complex reframing group $\Map(V(\Ga),\SL(2,\C))$ from
\eqref{cggdef}} and linear in the Higgs fields. In particular, the
dimension of the space of Higgs fields equals $3g-3$.

Now, following the Hitchin program (see \cite{H}) and sending
a Higgs field
\begin{equation}\label{hf}
   \phi\,:\;\wtd E_a\to\wtd E_a\otimes K_\Ga
\end{equation}
to the corresponding  quadratic differential\footnote{note that the constrains
\eqref{resmatcond} automatically imply the constrains
\eqref{rescond2}}
$\det\om(v)=-\om_{11}^2(v)-\om_{12}(v)\cdot\om_{21}(v)$, we get a
map
\begin{equation}\label{hinmap}
\pi\,:\;H^0(P_\Ga,\ad\wtd E_a\otimes K_\Ga)
 \to  H^0(P_\Ga,\Oh_{P_\Ga}(2\,K_\Ga))
\end{equation}

The restriction of the projective line bundle
$\PP(E)$ onto each $P_v$ is the quadric
$\PP(E)|_{P_v}=\PP_1\times P_v$, which admits two projections
$$\PP_1\stackrel{p_f}{\longleftarrow}\PP(E)|_{P_v}
       \stackrel{p_v}{\longrightarrow}P_v\,.
$$
We write $H=p_f^*\Oh_{\PP_1}(1)$ for the Grothendieck generator
such that $R^0{p_v}_*H=E|_{P_v}$. Then any Higgs field \eqref{hf}
is obtained as the direct image $R^0{p_v}_*$ of some homomorphism
$$\Phi\,:\;H\to H\otimes p_v^*\Oh_{P_v}(1)\,.
$$
Since $\phi$ is traceless, its restriction $\phi_p:\PP(E_p)\to\PP(E_p)$
onto a fiber over $p\in P_v$ is either a linear automorphism
$\PP(E_p)\stackrel\sim\to\PP(E_p)$, which has two distinct fixed
points $p_1,p_2\in\PP(E_p)$, or a degenerated map $\phi_p$, which
contracts $\PP(E_p)$ to one point $p_0\in\PP(E_p)$. Thus, we get a
ramified double covering
\begin{equation}\label{dc}
  w\,:\;\wtd P_v \to P_v
\end{equation}
which depends on and encodes the initial Higgs field $\phi$.

Let us suppose for simplicity that for each component $P_v$, $v\in
V(\Ga)$, the quadratic differential $\det\om(v)$ has two distinct
zero points $z_1(v),z_2(v)\in
P_v\setminus\{p_{e_1},p_{e_1},p_{e_1}\}$, where
$\{e_1,e_2,e_3\}=S(v)$. Then the ramification divisor of the double
covering \eqref{dc} is
\begin{equation}\label{rd}
   (\det\om(v))_0=z_1(v)+z_2(v)\subset P_v
\end{equation}
and the intersection number $\wtd P_v\cdot p_v^{-1}(p_e)=2$
on the quadric $\PP_1\times P_v$.

\begin{dfn}
  The pair $(\wtd P_v,L_\phi)$, where $L_{\phi}=H|_{\wtd
  P_v}\in\Pic\wtd P_v$, is called a {\it spectral data\/} of the
  Higgs field \eqref{hf}.
\end{dfn}

The Higgs field $\phi$ is reconstructed
from its spectral data as follows.
Twisting the standard adjoint sequence of sheaves
$$0\to\Oh(-\wtd P_v)\to\Oh\to\Oh_{\wtd P_v}\to0
$$
on the quadric $\PP_1\times P_v$ by the Grothendieck line bundle
$H$, we get an exact triple
\begin{equation}\label{etr}
0\to\Oh(-\wtd P_v)\otimes H\to H\to\Oh_{\wtd P_v}(H)\to0\;.
\end{equation}
Restricting the first term to the fiber, we get
$$\Oh(-\wtd P_v)(H)|_{\PP(E_p)}=\Oh_{\PP_1}(-2)\otimes\Oh_{\PP_1}(1)
                               =\Oh_{\PP_1}(-1)\,.
$$
Hence, applying $R{f_v}_*$ to \eqref{etr}, we get an isomorphism
$$0\to E\to R^0{f_v}_*(\Oh_{\wtd P_v}(H))\to0\,
$$
which coincides with \eqref{hf}.

Now let us switch on the gluing procedure prescribed by the
framing data $a:\vec E(\Ga)\to\SL(2,\C)$ to combine the spectral
curves $\wtd P_v$ together. Consider each
$\wtd P_v$ as projective line equipped with an
involution $i_v:\wtd P_v\to\wtd P_v$ such that
$\wtd P_v\stackrel w\to P_v=\wtd P_v/i_v$ is the double covering
\eqref{dc}. Then the ramification points \eqref{rd}
become the fixed points of $i_v$ and 3 picked points
$p_{e_i}\in P_v$ turn to 6 points $p_{e_i}^\pm\in\wtd P_v$
forming a triple of $i_v$-conjugated pairs. Gluing $\wtd P_{v_1}$
with $\wtd P_{v_2}$ in 2 points $p_e^{\pm}$ for each
$e=(v_1,v_2)\in S(v_1)\cap S(v_2)\subset E(\Ga)$ and
all $v_1,v_2\in V(\Ga)$, we get a reducible curve
$\wtd P_\Ga(\phi)$ of arithmetical genus  $4g-3$ with
the involution
$$i_\phi\,:\;\wtd P_\Ga(\phi)\to\wtd P_\Ga(\phi)
$$
such that the quotient $\wtd P_\Ga(\phi)/i_\phi=P_\Ga$ is
the original ll-curve we have started with. Under this procedure,
the local line bundles $L_\phi$ on $\wtd P_v$ are glued
to the global line bundle
$$L_\phi\in\Pic\wtd P_\Ga(\phi)\;,
$$
which restricts onto each $\wtd P_v$ as $L_\phi|_{\wtd
P_v}=\Oh_{\wtd P_v}(1)$.

Now, using the geometric interpretation of the Hitchin systems via
Prym varieties developed in \cite{T4}, \cite{T5} one can easily
prove the following

\begin{prop}
  The initial Higgs field \eqref{hf} on $P_\Ga$ is uniquely
  recovered from the triple $(\wtd P_\Ga(\phi),i_\phi,L_\phi)$;

and the fiber of the Hitchin map \eqref{hinmap}:
  $$\pi\,:\;H^0(P_\Ga,\ad\wtd E_a\otimes K_\Ga)
    \to  H^0(P_\Ga,\Oh_{P_\Ga}(2\,K_\Ga))
  $$
  over a regular quadratic differential
  $\Om\in H^0(P_\Ga,\Oh_{P_\Ga}(2 K_\Ga))$ coincides with the
  Prym variety:
  $$\pi^{-1}(\Om)=\mathrm{Prym}_w=
    \Pic_{1,...,1}\Bigl(\wtd P_\Ga(\phi)/w^*(\Pic(P_\Ga)\Bigr)\;.
  $$
\end{prop}
$$
$$
$$
$$
$$
$$
{\it Andrey Tyurin died in Bonn at 27th of October
going by bus to the Institute. The death was
instantaneous and unexpected. The text looked  
unfinished but it was
 decided that really it contains the result
and the argument to prove it; one just needs to
remember two classical papers listed below and
combine these ones with [T1]. Thus it should be 
strongly recommended to use these sources for understanding 
of the present text, just slightly edited by Alexei Gorodentsev
and Nikolai Tyurin.} 

$$
$$
$$
$$

\end{document}